\documentclass[12pt,reqno]{amsart}

\usepackage[latin1]{inputenc}
\usepackage{amsmath}
\usepackage{amsfonts}
\usepackage{amssymb}
\usepackage{graphics}
\usepackage{enumerate}
\usepackage{amssymb,amsmath,amsthm,amscd,latexsym,verbatim,graphicx,amsfonts}

\usepackage{amscd}
\usepackage{amsmath}
\usepackage{amssymb}
\usepackage{amsthm}
\usepackage{latexsym}
\usepackage{verbatim}

\theoremstyle{plain}
\newtheorem{theorem}{Theorem}[section]

\theoremstyle{definition}

\newtheorem{pr}[theorem]{Problem}
\theoremstyle{remark}

\newcommand{\kri}{k\rightarrow\infty}
\newcommand{\tri}{t\rightarrow\infty}

\newcommand{\bbC}{\mathbb{C}}

\newcommand{\bbD}{\mathbb{D}}

\newcommand{\bbN}{\mathbb{N}}

\newcommand{\mco}{\mathcal{O}}

\newcommand{\mct}{\mathcal{T}}

\title[]{Hyponormal Toeplitz Operators on the Bergman Space of the Disk}
\author[]{Nicole Revilla $\&$ Brian Simanek}
\date{}

\begin{document}
\maketitle

\begin{abstract}
We consider Toeplitz operators with bounded symbol acting on the Bergman space of the unit disk and assess their hyponormality.  We will mainly be concerned with the symbol $\varphi(z)=z^{n}|z|^{2s}+a(t)\bar{z}^{m}|z|^{2t}$, where $s$ and $t$ are positive real numbers and $m$ and $n$ are natural numbers.  The main goal is to understand how large $|a(t)|$ can be for this operator to be hyponormal and we will answer this question for large values of $t$.  We also correct a typo from a 2019 paper of Fleeman and Liaw concerning the norm of the commutator of the Toeplitz operator with symbol $z^m\bar{z}^n$ when $m>n$.
\end{abstract}

\vspace{4mm}

\footnotesize\noindent\textbf{Keywords:} Hyponormal operator, Toeplitz Operator, Bergman Space

\vspace{2mm}

\noindent\textbf{Mathematics Subject Classification:} Primary 47B20; Secondary 47B35

\vspace{2mm}

\normalsize

\section{Introduction}\label{Intro}


Consider the normalized area measure $dA$ on the unit disk and the corresponding Bergman space, which we denote by $A^2(\bbD)$.  Recall that this space is defined as follows
\[
A^2(\bbD)=\left\{f:\int_{\bbD}|f|^2dA<\infty,\, f\mathrm{\, is\, analytic\, in\, }\bbD\right\}=\left\{f(z)=\sum_{n=0}^{\infty}a_nz^n:\sum_{n=0}^{\infty}\frac{|a_n|^2}{n+1}<\infty\right\}
\]
and is equipped with the inner product
\[
\left\langle f,g\right\rangle=\int_{\bbD}f\bar{g}\,dA,\qquad\mbox{or}\qquad \left\langle\sum_{n=0}^{\infty}a_nz^n,\sum_{n=0}^{\infty}b_nz^n\right\rangle=\sum_{n=0}^{\infty}\frac{a_n\bar{b}_n}{n+1}.
\]

A bounded operator $T$ is said to be \textit{hyponormal} if $[T^*,T]\geq0$, where $T^*$ denotes the adjoint of $T$.  
The operators we will consider in this paper are Toeplitz operators whose symbol is in the space $L^{\infty}(\bbD)$.  If $\varphi\in L^{\infty}(\bbD)$, then the operator $T_{\varphi}:A^2(\bbD)\rightarrow A^2(\bbD)$ is defined by
\[
T_{\varphi}(f)=P(\varphi f),
\]
where $P:L^2(\bbD,dA)\rightarrow A^2(\bbD)$ denotes the orthogonal projection.  We are interested in understanding what symbols $\varphi$ yield Toeplitz operators $T_{\varphi}$ that are hyponormal.

Investigations of this type have a lengthy history and many known results have a convenient perturbative interpretation.  Specifically, many known results address the question: if $T_{\varphi}$ is hyponormal and $\psi\in L^{\infty}(\bbD)$, then for what values of $c\in\bbC$ is it true that $T_{\varphi+c\psi}$ is also hyponormal?  For many functions $\varphi$ and $\psi$, it is true that the answer to this question depends only on $|c|$ and states that $|c|$ can be no bigger than some particular value (that depends on $\varphi$ and $\psi$).  Results of this kind can be found in \cite{CC,FL,KL21,LS,Sadraoui,Sadraoui21,SimHypo}.

Our results are improvements on some results presented in \cite{FL}.  It is known that $T_{z^m\bar{z}^n}$ is hyponormal if and only if $m\geq n$ (see \cite[Theorem 2.3]{KL21}).  Therefore, we can ask: if $\varphi=z^{p+q}\bar{z}^p+a\bar{z}^{j+k}z^j$, how large can $|a|$ be and still have $T_{\varphi}$ be hyponormal?  The following result, which provides a partial answer to this question, is taken directly from \cite{FL}.

\begin{theorem}\label{FL10}\cite[Theorem 10]{FL}
Fix $\delta\in\bbN$.  For every $n\in\bbN$, there exists $j\in\bbN$ so that $T_{\varphi}$ with symbol $\varphi(z)=z^{n+\delta}\bar{z}^n+\frac{1}{2j+\delta}\bar{z}^{j+\delta}z^j$ is hyponormal.
\end{theorem}

The wording of this theorem is suboptimal in the sense that something stronger is known.  Indeed, the following result was proven in \cite{SimHypo} (see also \cite[Theorem 2.4]{KL21}).

\begin{theorem}\label{shcor2}\cite[Corollary 2]{SimHypo}
If $n\in\bbN$ and $s\in[0,\infty)$, then $T_{|z|^s(z^n+a\bar{z}^n)}$ is hyponormal if and only if $|a|\leq1$.
\end{theorem}

Thus, in Theorem \ref{FL10} one can set $j=n$ and then replace $(2j+\delta)^{-1}$ by any complex constant in the closed unit disk and the corresponding operator $T_{\varphi}$ will be hyponormal.  However, the proof from \cite{FL} actually produces the following result:

\begin{theorem}\label{FL10+}
Fix $\delta\in\bbN$.  If $n\in\bbN$ is fixed, then for all $j\in\bbN$ sufficiently large, the operator $T_{\varphi}$ with symbol $\varphi(z)=z^{n+\delta}\bar{z}^n+\frac{1}{2j+\delta}\bar{z}^{j+\delta}z^j$ is hyponormal.
\end{theorem}

\noindent\textit{Remark.}  A related result was proven in \cite[Theorem 2.5]{KL21}.

\medskip

This wording puts a different spin on the problem by asking how the constant $j$ influences the magnitude of the perturbations that preserve hyponormality for these symbols and specifically addresses this behavior as $j\rightarrow\infty$.  What Theorem \ref{FL10+} does not address is the influence of $n$ on this problem.  Our first theorem improves Theorems \ref{FL10} and \ref{FL10+} by considering a more general range of parameters, refining the constants involved, and showing that our result is sharp.

\begin{theorem}\label{NewFL10}
Fix $s\in(0,\infty)$ and $m,n\in\bbN$.  Suppose $a(t):(0,\infty)\rightarrow\bbC$ and consider the operator $T_{\varphi}$ with symbol $\varphi(z)=z^{n}|z|^{2s}+a(t)\bar{z}^{m}|z|^{2t}$.


\noindent(i)  If
\[
\limsup_{\tri}|ta(t)|< \frac{n+2s}{2},
\]
then for all sufficiently large $t$, it holds that $T_{\varphi}$ is hyponormal.


\noindent(ii) If
\[
\liminf_{\tri}|ta(t)|> \frac{n+2s}{2},
\]
then for all sufficiently large $t$, it holds that $T_{\varphi}$ is not hyponormal.
\end{theorem}

In the notation of Theorem \ref{FL10}, we see that for large $j$, one could replace $(2j+\delta)^{-1}$ by $C/(2j)$, where $C$ satisfies $|C|<\delta+2n$.  Notice that in Theorem \ref{NewFL10}, it is clear how the parameters $n$, $s$, and $t$ influence which values of $a(t)$ yield hyponormal operators, but it is not clear how $m$ influences these values.  What is likely true is that $m$ will influence exactly how large is ``sufficiently large," i.e. for any specific choice of $n$, $m$, and $s$, all three parameters will determine precisely which values of $t$ do or do not yield hyponormal operators.

It follows easily from our construction that there are examples of functions $a(t)$ for which
\[
\lim_{\tri}ta(t)= \frac{n+2s}{2}
\]
and the operator $T_{\varphi}$ is hyponormal for some large values of $t$ and not hyponormal for some other large values of $t$.  Indeed, such a function $a(t)$ could appropriately alternate between $(1+\epsilon_{\ell})(n+2s)/(2t)$ and $(1-\epsilon_{\ell})(n+2s)/(2t)$ for a positive sequence $\{\epsilon_{\ell}\}_{n\in\bbN}$ converging monotonically to zero as $\ell\rightarrow\infty$.

At first glance, one might find the result of Theorem \ref{NewFL10} surprising.  Indeed, the symbol $z^n|z|^{2s}$ yields a hyponormal operator while the symbol $\bar{z}^m|z|^{2t}$ does not.  Notice that as $\tri$, the function $\bar{z}^m|z|^{2t}$ tends to zero on $\bbD$.  Therefore, if we think of this function as being the perturbation term in the symbol of Theorem \ref{NewFL10}, then as $\tri$ the perturbation is becoming smaller and hence one might expect $|a(t)|$ can grow as $\tri$ and still preserve hyponormality.  However, several results in the literature suggest that it is really the boundary behavior of the symbol (or its derivative) that is most relevant in determining hyponormality (see for example \cite[Proposition 1.4.3]{Sadraoui}, \cite[Theorem 2]{Sadraoui21}, and \cite[Cor. to Thm. 4]{AC}) .  Since $\bar{z}^m|z|^{2t}$ does not go to zero on the boundary of the unit disk as $\tri$, the result of Theorem \ref{NewFL10} becomes less surprising.

In light of Theorem \ref{NewFL10}, it is natural to ask the complementary question of what happens as $t\rightarrow0$.  In this case, the values of $|a(t)|$ that yield hyponormal operators do not blow up as $t^{-1}$, but instead remain bounded.  It will be a simple matter to prove the following theorem.

\begin{theorem}\label{tzero}
Fix $s\in(0,\infty)$ and $m,n\in\bbN$ and consider the operator $T_{\varphi}$ with symbol $\varphi(z)=z^{n}|z|^{2s}+a\bar{z}^{m}$. If $T_{\varphi}$ is hyponormal, then
\[
|a|^2\leq\min\left\{\frac{(m+1)(n+1)}{(n+s+1)^2},\frac{n(n+2s)}{m^2}\right\}.
\]
\end{theorem}

\noindent\textit{Remark.}  If we let $s\rightarrow0$ in Theorem \ref{tzero}, then we recover the necessary conditions in \cite[Proposition 1.4.4]{Sadraoui}.

\medskip

We can also consider the setting of Theorem \ref{NewFL10} in which $s$ is related to $t$ and $n$ is related to $m$.  This naturally leads us to an improvement of \cite[Theorem 2.7]{KL21}.  If we fix $m,q\in\bbN$ such that $m-q-1\geq0$, then we can consider the operator $T_{\varphi}$ with symbol $\varphi(z)=z^m\bar{z}^{m-1}+a\bar{z}^{m-q}z^{m-q-1}$.  

\begin{theorem}\label{NewKL2.7}
Fix $m,q\in\bbN$ such that $m-q-1\geq0$.  If the operator $T_{\varphi}$ with symbol $\varphi(z)=z^m\bar{z}^{m-1}+a\bar{z}^{m-q}z^{m-q-1}$ is hyponormal, then
\[
|a|\leq\frac{m-q+1}{m+1}.
\]
\end{theorem}

\noindent\textit{Remark.}  The special case of Theorem \ref{NewKL2.7} in which $q=1$ was proven in \cite[Theorem 2.7]{KL21}.

\medskip

The second part of the paper will correct a mistake in \cite[Theorem 7]{FL}.  In that result, the authors considered the commutator $[T^*_{z^m\bar{z}^n},T_{z^m\bar{z}^n}]$ with $m>n$ (when $m=n$, the commutator is $0$).  The operator $T_{z^m\bar{z}^n}$ is hyponormal and the proof in \cite{FL} correctly identifies the eigenvalues of the commutator as $\{\lambda_k\}_{k=0}^{\infty}$, where
\[
\lambda_k=\begin{cases}
(k+1)\left(\frac{k+m-n+1}{(k+m+1)^2}\right) \qquad\qquad\qquad\qquad\qquad & \mbox{if}\quad 0\leq k<m-n\\
(k+1)\left(\frac{k+m-n+1}{(k+m+1)^2}-\frac{k+n-m+1}{(k+n+1)^2}\right) & \mbox{if}\quad k\geq m-n.
\end{cases}
\]
The norm of the commutator is then the maximum of these eigenvalues, but the proof in \cite{FL} incorrectly asserts that $(k+1)\left(\frac{k+m-n+1}{(k+m+1)^2}-\frac{k+n-m+1}{(k+n+1)^2}\right)$ is a ``monotonically decreasing function in $k$" in the domain $[m-n,\infty)$.  Unfortunately, this leads to an incorrect identification of which eigenvalue is largest and hence yields an incorrect value for the norm of the commutator.

Identifying the largest eigenvalue is a difficult task, but we can prove the following result, which conveniently identifies pairs $(m,n)$ with $m\geq n$ for which the assertion in \cite{FL} is true and describes when the function $(x+1)\left(\frac{x+m-n+1}{(x+m+1)^2}-\frac{x+n-m+1}{(x+n+1)^2}\right)$ is not a monotonically decreasing function of the continuous variable $x$ in $[m-n,\infty)$.

\begin{theorem}\label{NewFL7}
Suppose $m>n$.  The function
\[
F(x)=(x+1)\left(\frac{x+m-n+1}{(x+m+1)^2}-\frac{x+n-m+1}{(x+n+1)^2}\right)
\]
is a monotonically decreasing function of $x\in[m-n,\infty)$ if and only if
\begin{equation}\label{mncon}
(1+m)^3(m+n-mn+2m^2-n^2)<(2m+1-n)^3(m+n-mn+m^2).
\end{equation}
If $F$ is not monotonically decreasing on $[m-n,\infty)$, then it has a unique critical point in this interval, where $F$ attains a local maximum.
\end{theorem}

The condition \eqref{mncon} is easy to verify for any particular choice of $m\geq n$.  Indeed, we can easily produce Figure 1, where the shaded region shows those pairs $(m,n)$ up to $10000$ that do not satisfy this condition.  We see that the shaded region is bounded by two approximately straight lines.  The top one is $m=n$ and we will prove that the bottom boundary curve is approximately a line with equation $m/n\approx1.61\ldots$.

\begin{figure}
  \centering
  \includegraphics[width=0.63\textwidth]{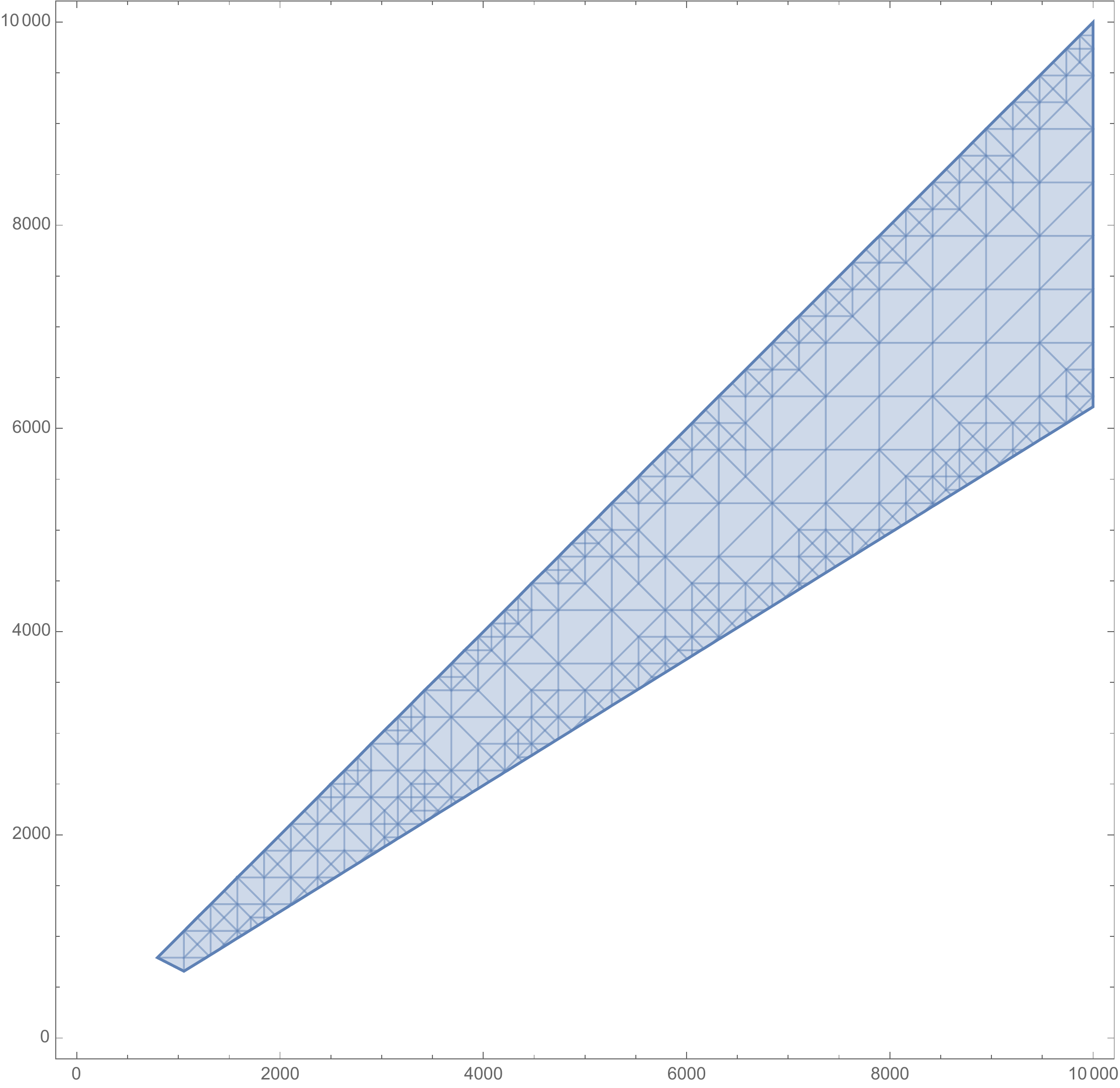}
\caption{The shaded region shows those pairs $(m,n)$ up to $10000$ that do not satisfy condition \eqref{mncon}. The horizontal axis is $m$ and the vertical axis is $n$.}
\label{fig:1}       
\end{figure}

Strictly speaking, the condition that $F(x)$ from Theorem \ref{NewFL7} has a local maximum in $(m-n,\infty)$ does not imply that the sequence of eigenvalues $\{\lambda_k\}_{k=m-n}^{\infty}$ is not a decreasing sequence.  If $x^*$ is a critical point of $F$ in $[m-n,\infty)$, then the maximum of $\{\lambda_k\}_{k=m-n}^{\infty}$ is attained for $k=\lfloor x^*\rfloor$ or $k=\lceil x^*\rceil$.  It could happen that the $x^*\in(m-n,m-n+1)$ and the maximum of $\{\lambda_k\}_{k=m-n}^{\infty}$ is attained when $k=m-n$ (this happens, for example, when $m=5$ and $n=4$).  However, it is easy to find examples of pairs $(m,n)$ with $m>n$ so that the maximum of $\{\lambda_k\}_{k=m-n}^{\infty}$ is attained at a value of $k$ different from $m-n$.  For example, if $m=8$ and $n=7$, then $\max_{k\in\bbN_0}\{\lambda_k\}=\lambda_3$.

In light of the error in \cite[Theorem 7]{FL}, we must also address \cite[Corollary 2]{FL}, which states
\begin{equation}\label{NewCor2}
\max_{k}\{\lambda_k\}\leq\frac{1}{2}.
\end{equation}
While the calculations that justify this claim in \cite{FL} are based on the aforementioned error, we will show that this estimate is still correct and in fact one can make the inequality in \eqref{NewCor2} strict.

\medskip

Now that we have introduced our results and notation, we can turn our attention to the proofs.  In the next section, we will prove Theorem \ref{NewFL10}.  In subsequent sections we will prove Theorems \ref{tzero}, \ref{NewKL2.7}, and \ref{NewFL7}.  

\bigskip

\noindent\textbf{Acknowledgements.} The second author gratefully acknowledges support from the Simons Foundation through collaboration grant 707882.

\section{Proof of Theorem \ref{NewFL10}}\label{NewFL10proof}

Recall from \cite[Section 3]{SimHypo} that we may define (we change the notation from $\sigma_k'$ to $\sigma_k$, etc. for convenience)
\begin{align*}
\sigma_k&=
\begin{cases}
\frac{k+n+1}{(n+k+s+1)^2}\qquad\qquad\qquad\qquad\qquad\qquad\qquad & 0\leq k< n\\
\frac{k+n+1}{(n+k+s+1)^2}-\frac{k-n+1}{(k+s+1)^2} & k\geq n
\end{cases}
\\
\omega_k&=
\begin{cases}
-\frac{k+m+1}{(m+k+t+1)^2}\qquad\qquad\qquad\qquad\quad\qquad\qquad & 0\leq k< m\\
\frac{k-m+1}{(k+t+1)^2}-\frac{k+m+1}{(m+k+t+1)^2} & k\geq m
\end{cases}
\\
\delta_k&=
\frac{k+n+1}{(k+n+s+1)(k+n+m+t+1)}-\frac{k+m+1}{(k+n+m+s+1)(k+m+t+1)}
\end{align*}

It was shown in \cite{SimHypo} that $T_{\varphi}$ is hyponormal if and only if
\begin{equation}\label{hncon}
|a(t)|^2\sum_{k=0}^{\infty}\omega_k(t)|u_k|^2-2|a(t)|\sum_{k=0}^{\infty}|\delta_k(t)u_k\bar{u}_{k+n+m}|+\sum_{k=0}^{\infty}\sigma_k|u_k|^2>0
\end{equation}
for all $\{u_k\}_{k=0}^{\infty}$ satisfying $\sum_{k=0}^{\infty}|u_k|^2/(k+1)<\infty$.  It is also known that each $\sigma_k>0$ and each $\omega_k<0$.

Since we will be sending $\tri$, let us think of $\omega_k$ and $\delta_k$ as functions of $t$.  When we write $x_k(t)=o(1)$, we mean that if $\epsilon>0$ is given, then there is a constant $A_{\epsilon}>0$ so that if $t,k>A_\epsilon$, then $|x_k(t)|<\epsilon$.  Observe that as $\kri$
\begin{align}\label{sigmaform}
\nonumber\sigma_k&=\frac{n(n^2+n(2s+k+1)+2s(k+s+1))}{(k+s+1)^2(k+n+s+1)^2}\\
&=\frac{n(n+2s)}{k^3}\left(1+\frac{(n^2-2s^2)/(n+2s)-2n-2s-3}{k}+O(k^{-2})\right)
\end{align}
We also observe that
\begin{align}\label{omegaform}
\nonumber\omega_k(t)&=\frac{-m(m^2+m(2t+k+1)+2t(k+t+1))}{(k+t+1)^2(k+m+t+1)^2}\\
&=\frac{-2mt}{(k+t)^3}\left(1+\frac{m(k+m+1)/2t-m-3}{k+t}+\mco\left(\frac{1+k/t}{(k+t)^2}\right)\right)
\end{align}
and
\begin{align}\label{deltaform}
\nonumber\delta_k(t)&=\frac{nt(k+n+s+1)-m^2s-ms(k+t+1)}{(k+n+s+1)(k+n+m+t+1)(k+n+m+s+1)(k+m+t+1)}\\
&=\frac{nt}{k(k+t)^2}\left(1-\frac{1+n+m+s+\frac{ms}{n}}{k}-\frac{2m+n+2}{k+t}-\frac{ms/n}{t}+\mco\left(\frac{1}{k^2}+\frac{1}{kt}\right)\right)
\end{align}
as $\tri$ and the error estimates are uniform in $k\in\bbN$.

\subsection{Proof of Part i}

From \eqref{omegaform} and \eqref{sigmaform}, we see that as $\tri$
\begin{align}\label{omegat}
\nonumber\left|\frac{\omega_k(t)}{t^2}\right|&\leq \frac{2m}{t(k+t)^3}\left(1+\mco(t^{-1})\right)\leq\frac{2m(1+o(1))}{n(n+2s)t}\sigma_k
\end{align}

We can also use the above formulas to see that if $t>s$, then as $\kri$
\begin{equation}\label{deltat}
\left|\frac{\delta_k(t)}{t}\right|\leq\frac{n(1+o(1))}{k(k+t)^2}\leq\frac{(1+o(1))}{(n+2s)}\sigma_k
\end{equation}

It is clear that for each $k\in\bbN_0$, it holds that $\delta_k(t)\rightarrow0$ and $\omega_k(t)\rightarrow0$ as $\tri$ and that $\sigma_k/\sigma_{k+1}\rightarrow1$ as $\kri$.  Thus, given any $\epsilon>0$, we may find $M>0$ and $\mct>0$ so that
\begin{align*}
|\delta_k(t)|+|\omega_k(t)|<\min\{\sigma_k,\sigma_{k+n+m}\}\qquad\qquad &k<M\\
\left|\frac{\delta_k(t)}{t}\right|\leq\frac{(1+\epsilon)}{(n+2s)}\min\{\sigma_k,\sigma_{k+n+m}\},\qquad\qquad &k\geq M\\
\left|\frac{\omega_k(t)}{t^2}\right|\leq\frac{(2m+\epsilon)}{n(n+2s)t}\sigma_k,\qquad\qquad &k\geq M
\end{align*}
for all $t>\mct$.

Assume from now on that $t>\mct$.  Applying the Cauchy-Schwartz inequality shows
\begin{align}\label{firstbelow}
\nonumber&|a(t)|^2\sum_{k=0}^{\infty}\omega_k(t)|u_k|^2-2|a(t)|\sum_{k=0}^{\infty}|\delta_k(t)u_k\bar{u}_{k+n+m}|+\sum_{k=0}^{\infty}\sigma_k|u_k|^2\\
\nonumber&\qquad>\sum_{k=0}^{\infty}\sigma_k|u_k|^2\bigg(1-2|a(t)|\frac{\sqrt{\sum_{k=M}^{\infty}|\delta_k(t)||u_{k}|^2}\sqrt{\sum_{k=M}^{\infty}|\delta_k(t)||u_{k+n+m}|^2}}{\sum_{k=0}^{\infty}\sigma_k|u_k|^2}\\
\nonumber&\qquad\qquad\qquad\qquad\qquad\qquad\qquad\qquad\qquad-|a(t)|^2\frac{\sum_{k=M}^{\infty}|\omega_k(t)||u_k|^2}{\sum_{k=0}^{\infty}\sigma_k|u_k|^2}-2|a(t)|-|a(t)|^2\bigg)\\
&\qquad\geq\sum_{k=0}^{\infty}\sigma_k|u_k|^2\bigg(1-2|a(t)|\frac{(1+\epsilon)t}{n+2s}-|a(t)|^2t\frac{2m+\epsilon}{n(n+2s)}-2|a(t)|-|a(t)|^2\bigg)
\end{align}
Our assumptions imply that for all large $t$, there is $\kappa\in(0,1)$ so that $|ta(t)|\leq\frac{1}{2}\kappa(n+2s)$.  Therefore, we may use our conditions on $M$ to bound the right-hand side of \eqref{firstbelow} from below by
\begin{align}\label{secondbelow}
\sum_{k=0}^{\infty}\sigma_k|u_k|^2\bigg(1-\kappa(1+\epsilon)-\mco(t^{-1})\bigg)
\end{align}
as $\tri$.  If we choose $\epsilon$ small enough so that $\kappa(1+\epsilon)<1$, then \eqref{secondbelow} is strictly positive for all sufficiently large $t$ and this estimate is independent of the sequence $\{u_k\}_{k=0}^{\infty}$.  Therefore, the operator $T_{\varphi}$ is hyponormal.

\subsection{Proof of Part ii}

Now suppose $|a(t)|>\frac{\eta(n+2s)}{t}$ for some $\eta>1/2$ when $t$ is large.  Fix some $k_1,k_2\in\bbN$ and consider the sequence
\[
u_k=\begin{cases}
1\qquad\qquad&\mbox{if}\quad k_1\leq k\leq k_1+k_2\\
0 & \mbox{otherwise}.
\end{cases}
\]

Pick $\epsilon>0$ small enough and $k_2$ large enough so that
\begin{equation}\label{epsdef}
2\eta\frac{(k_2-n-m+1)(1-\epsilon)}{(k_2+1)(1+\epsilon)}>1.
\end{equation}
We will show that \eqref{hncon} is violated by this choice of $\{u_n\}_{n=0}^{\infty}$ for all sufficiently large choices of $k_1$ and $t=\sqrt{k_1}$.

A short calculation reveals that
\[
\delta_k(t)=\frac{t(n^2-ms+n(s+k+1))-ms(m+k+1)}{P_k},
\]
where $P_k$ is positive for all choices of $k,m,n,s$, and $t$.  Thus, if $k$ and $t$ are (independently) large enough, then it holds that $\delta_k(t)>0$.

If $k_1$ is large enough and we set $t=\sqrt{k_1}$, then the left-hand side of \eqref{firstbelow} is bounded above by
\[
\sum_{k=k_1}^{k_1+k_2}\sigma_{k}-\frac{2\eta(n+2s)}{\sqrt{k_1}}\sum_{k=k_1}^{k_1+k_2-n-m}|\delta_{k}\left(\sqrt{k_1}\right)|-\frac{\eta^2(n+2s)^2}{k_1}\sum_{k=k_1}^{k_1+k_2}|\omega_{k}\left(\sqrt{k_1}\right)|
\]
Assume $k_1$ is large enough to ensure that $\delta_k(t)>0$ for all $t>\sqrt{k_1}$ and $k>k_1$.  We can upper-bound the above expression by
\begin{align}\label{maxminii}
&(k_2+1)\max_{k_1\leq j\leq k_1+k_2}\sigma_{j}-\frac{2\eta(n+2s)}{\sqrt{k_1}}(k_2-n-m+1)\min_{k_1\leq j\leq k_1+k_2-n-m}\delta_{j}\left(\sqrt{k_1}\right)
\end{align}

By using the estimate \eqref{deltaform}, we see that if $k_1$ is large enough, then \eqref{maxminii} is bounded from above by
\begin{equation}\label{lastii}
(k_2+1)\frac{n(n+2s)(1+\epsilon)}{k_1^3}-\frac{2\eta(n+2s)}{\sqrt{k_1}}\cdot\frac{(k_2-n-m+1)n\sqrt{k_1}(1-\epsilon)}{(k_1+k_2-n-m)(k_1+\sqrt{k_1}+k_2-n-m)^2}
\end{equation}
This simplifies to
\[
\frac{(k_2+1)(1+\epsilon)n(n+2s)}{k_1^3}\left[1-\frac{2\eta(k_2-n-m+1)(1-\epsilon)k_1^3}{(k_2+1)(1+\epsilon)(k_1+k_2-n-m)(k_1+\sqrt{k_1}+k_2-n-m)^2}\right]
\]
Notice that
\[
\frac{2\eta(k_2-n-m+1)(1-\epsilon)k_1^3}{(k_2+1)(1+\epsilon)(k_1+k_2-n-m)(k_1+\sqrt{k_1}+k_2-n-m)^2}=2\eta\frac{(k_2-n-m+1)(1-\epsilon)}{(k_2+1)(1+\epsilon)}+o(1)
\]
as $k_1\rightarrow\infty$ and so this expression is strictly larger than $1$ when $k_1$ is large, by \eqref{epsdef}.  Therefore, if $k_1$ is sufficiently large, then the expression \eqref{lastii} is negative and hence the operator $T_{\varphi}$ is not hyponormal with this choice of $a(t)$.


\section{Proof of Theorem \ref{tzero}}\label{tzeroproof}

We again use \eqref{hncon} and the quantities $\sigma_k$, $\delta_k$, and $\omega_k$ but with $t=0$.  If we plug in the test vector $e_k=(0,0,\ldots,0,1,0,\ldots)^T$ with a $1$ in position $k$ and zeros elsewhere, then \eqref{hncon} implies that hyponormality of $T_{\varphi}$ requires
\[
|a|^2\leq\frac{\sigma_k}{-\omega_k}.
\]
If we set $k=0$, then the right-hand side of this inequality becomes $(m+1)(n+1)/(n+s+1)^2$.  If we send $\kri$, then the right-hand side is $n(n+2s)m^{-2}(1+O(k^{-1}))$.  Since the inequality must hold for all $k$, the desired result follows.

\section{Proof of Theorem \ref{NewKL2.7}}\label{NewKL2.7proof}

Recall the notation $\lambda_k$ from Section \ref{Intro} and write $\lambda_k=\lambda_k(m,n)$.  By reproducing the argument from the proof of \cite[Theorem 2.7]{KL21}, one finds that hyponormality of $T_{\varphi}$ implies
\[
|a|^2\leq\min\left\{\left(\frac{m-q+1}{m+1}\right)^2,\min_{k\in\bbN}\left\{\frac{\lambda_k(m,m-1)}{\lambda_k(m-q,m-q-1)}\right\}\right\}
\]
A straightforward (but lengthy) calculation reveals that $\lambda_k(m,m-1)/\lambda_k(m-q,m-q-1)$ is a decreasing function of $k\in\bbN$.  Therefore, we conclude that
\[
|a|^2\leq\min\left\{\left(\frac{m-q+1}{m+1}\right)^2,\frac{\frac{3}{(m+2)^2}-\frac{1}{(m+1)^2}}{\frac{3}{(m-q+2)^2}-\frac{1}{(m-q+1)^2}}\right\}
\]
One can check that this minimum always equals $(m-q+1)^2/(m+1)^2$ when $m\geq q+1$ and the result follows.

\section{Proof of Theorem \ref{NewFL7}}\label{NewFL7proof}

We begin by calculating
\[
F'(x)=\frac{(x+n+1)^3((m+n)(x+m+1)-2mn)-(x+m+1)^3((m+n)(x+n+1)-2mn)}{(x+m+1)^3(x+n+1)^3}
\]
We wish to identify the zeros of this function in the interval $[m-n,\infty)$.  Thus, we will write
\[
F'(x+m-n)=\frac{P(x)}{Q(x)}.
\]
so that a positive zero of $P(x)$ corresponds to a critical point of $F$ in $(m-n,\infty)$.  It is straightforward to verify that $Q(x)>0$ when $x\geq0$ and
\[
P(x)=ax^3+bx^2+cx+d,
\]
where
\begin{align*}
a&=2(n^2-m^2)\\
b&=3(n-m)(3m^2-n^2+2m+2n)\\
c&=(n-m)(13m^3+18m^2+6m-13m^2n+6n-mn^2-6n^2+n^3)\\
d&=(1+m)^3(m+n-mn+2m^2-n^2)-(2m+1-n)^3(m+n-mn+m^2).
\end{align*}
It is clear from this presentation that $a<0$ and $b<0$ when $m>n$.  Also notice that
\[
c=(n-m)((13m^2-n^2)(m-n)+6(m+n)(m-n+1)+12m^2),
\]
which is clearly negative when $m>n$.  Thus, only $d$ can be either positive or negative and Descarte's Rule of Signs tells us that the number of positive critical points of $F$ will be $1$ if $d>0$ and $0$ if $d<0$.  Thus, the absence of a critical point is equivalent to $d<0$ and this is precisely the statement of Theorem \ref{NewFL7}.  The fact that the critical point (if it occurs) is a local maximum of $F$ follows from the fact that $P(0)>0$ when $d>0$, $Q(x)>0$ when $x\geq0$, and $P(x)\rightarrow-\infty$ as $x\rightarrow\infty$.

\bigskip

To approximate the curve that is the lower boundary of the shaded region in Figure 1, we look at the formula for $d$ and replace $m$ by $\alpha n$ (with $\alpha>1$).  Then
\[
d=(1+\alpha n)^3((\alpha+1)n+(2\alpha^2-\alpha-1)n^2)-((2\alpha-1)n+1)^3((\alpha+1)n+(\alpha^2-\alpha)n^2)
\]
This is a polynomial of degree $5$ in $n$ and if it is equal to $0$, then the leading coefficient must be zero.  The leading coefficient is
\[
\alpha(\alpha-1)(\alpha^2(2\alpha+1)-(2\alpha-1)^3).
\]
The values $\alpha=0$ and $\alpha=1$ are not in our parameter space and the only real solution of $\alpha^2(2\alpha+1)-(2\alpha-1)^3=0$ is $\alpha\approx1.61\ldots$, which proves our claim.

\bigskip

To justify the claim \eqref{NewCor2}, we must show that each $\lambda_k$ is at most $1/2$.  The calculations in \cite{FL} are valid for $0\leq k<m-n$ (and yield a strict inequality in \eqref{NewCor2}), so we will only consider $k\geq m-n$.  After clearing denominators, this reduces to showing that
\[
2(k+1)\left((k+m-n+1)(k+n+1)^2-(k+n-m+1)(k+m+1)^2\right)\leq(k+m+1)^2(k+n+1)^2
\]
when $k\geq m-n$.  Define
\begin{align*}
&R(x):=(x+m+1)^2(x+n+1)^2\\
&\qquad\qquad\qquad-\left[2(x+1)\left((x+m-n+1)(x+n+1)^2-(x+n-m+1)(x+m+1)^2\right)\right]
\end{align*}
so it suffices to show that $R(x)$ is positive on $[m-n,\infty)$.  We can write
\[
R(x+m-n)=x^4+\alpha x^3+\beta x^2+\gamma x+\delta,
\]
where
\begin{align*}
\alpha&=4+6m-2n\\
\beta&=6(2m+1)+(8m+6)(m-n)+3(m^2+n^2)\\
\gamma&=2\left(2m^3+3(m^2+n^2)+6m+2+(m-n)(m^2-mn+n^2+8m+3)\right)\\
\delta&=(m+1)^2(2m-n+1)^2-2(m-n)(m+n-1)(2m^2+m+n)
\end{align*}
It is clear from these formulas that $m>n$ implies $\alpha>0$, $\beta>0$, and $\gamma>0$.  To show that $\delta>0$ also, write $\delta=\delta(m,n)$ and notice that $\delta(n,n)>0$ and
\[
\frac{\partial}{\partial m}\delta(m,n)=2\left(3m^2n+5m^2+2n^2+7m+3+(m-n)(3mn+4m+4)\right)>0
\]
when $m\geq n$.  It follows that $\delta(m,n)>0$ when $m>n$ as desired.  Notice that this reproves the assertion from the proof of \cite[Corollary 2]{FL} that $\lambda_{m-n}\leq1/2$.

Combining all of these calculations shows $R(x)>0$ when $x\in[m-n,\infty)$ and hence the inequality \eqref{NewCor2} with strict inequality follows.  Even though this inequality is strict, it is clear that one cannot obtain a tighter bound that $1/2$ that is valid for all pairs $(m,n)$ with $m>n$.

\section{Future Work}\label{future}

The ultimate goal of this line of research is to find a characterization of the $L^{\infty}(\bbD)$ symbols $\varphi$ that yield hyponormal operators $T_{\varphi}$ acting on the Bergman space.  The analogous characterization in the Hardy space $H^2(\bbD)$ was obtained by Cowen in \cite{Cowen}, thought the Bergman space version currently seems out of reach.  Our new results naturally lead to several problems that would be interesting to solve and may provide valuable insight into the more general problem.

We begin with the most natural extension of Theorem \ref{NewFL10}.

\begin{pr}\label{prFL10}
Find all values of $n,m\in\bbN$ and $s,t\in(0,\infty)$ for which $T_{\varphi}:A^2(\bbD)\rightarrow A^2(\bbD)$ is hyponormal when
\[
\varphi(z)=z^n|z|^{2s}+\frac{n+2s}{2t+1}\bar{z}^m|z|^{2t}.
\]
\end{pr}

Based on our calculations, it is a difficult task to even conjecture a description of all $4$-tuples $(n,m,s,t)$ that solve Problem \ref{prFL10}.  Theorem \ref{tzero} suggests that the answer will depend on $m$ in a way that is not yet clear.





\medskip

Our next problem relates to Theorem \ref{NewFL7} and its solution would be a complete replacement for \cite[Theorem 7]{FL}.

\begin{pr}\label{prFL7}
If $m>n$, calculate the norm of the commutator $[T_{z^m\bar{z}^n}^*,T_{z^m\bar{z}^n}]$.
\end{pr}

We know that the norm is the largest of the eigenvalues $\{\lambda_k\}_{k=0}^{\infty}$, but a complete solution to the problem would describe - for each $k\in\{0,1,2,\ldots\}$ - the set of pairs $(m,n)$ with $m>n$ for which $\|[T_{z^m\bar{z}^n}^*,T_{z^m\bar{z}^n}]\|=\lambda_k$.  A helpful first step might be to describe those pairs $(m,n)$ for which $\|[T_{z^m\bar{z}^n}^*,T_{z^m\bar{z}^n}]\|=\lambda_{m-n-1}$.

\end{document}